\documentclass[11pt]{article}
\usepackage{latexsym,amsmath,amsfonts,graphics}
\begin{document}
\begin{center}
{\LARGE\textbf{Joint probability generating function for degrees of active/passive random intersection graphs}}\\
\bigskip
\bigskip
Yilun Shang\footnote{Department of Mathematics, Shanghai Jiao Tong University, Shanghai 200240, CHINA. email: \texttt{shyl@sjtu.edu.cn}}\\
\end{center}

\begin{abstract}
Correlations of active and passive random intersection graphs are
studied in this letter. We present the joint probability generating
function for degrees of $G^{active}(n,m,p)$ and
$G^{passive}(n,m,p)$, which are generated by a random bipartite
graph $G^*(n,m,p)$ on $n+m$ vertices.

\bigskip

\end{abstract}

\bigskip
\normalsize

\noindent{\Large\textbf{1. Introduction}}
\smallskip

Consider a set $V$ with $n$ vertices and another set of objects $W$
with $m$ objects. Define a bipartite graph $G^*(n,m,p)$ with
independent vertex sets $V$ and $W$. Edges between $v\in V$ and
$w\in W$ exist independently with probability $p$. The active random
intersection graph $G^{active}(n,m,p)$ derived from $G^*(n,m,p)$ is
defined on the vertex set $V$ with vertices $v_1,v_2\in V$ adjacent
if and only if there exists some $w\in W$ such that both $v_1$ and
$v_2$ are adjacent to $w$ in $G^*(n,m,p)$. Analogously, the passive
random intersection graph $G^{passive}(n,m,p)$ is defined on the
vertex set $W$ with vertices $w_1,w_2\in W$ adjacent if and only if
there is some $v\in V$ such that both $w_1$ and $w_2$ are adjacent
to $v$ in $G^*(n,m,p)$.

The models of random intersection graphs defined above were first
introduced in \cite{2,3}. Observe that the degree of a given vertex
$v\in V$ (or $w\in W$, respectively) in $G^*(n,m,p)$ is binomial
$Bin(m,p)$ (or $Bin(n,p)$, respectively) distributed, and by
interchanging the roles of $V$ and $W$, the active and passive
graphs are dual to each other essentially. The random intersection
graphs are generalized in \cite{5} by allowing an arbitrary
distribution for vertex degree in $G^*(n,m,p)$ instead of merely
binomial. Extended in such a way, the active and passive graphs
reveal different properties including degree distributions; see
\cite{5} for details. Intersection graphs are relationship graphs
and widely applied in various fields such as investigation of secure
wireless sensor networks \cite{8}, modeling of social networks
\cite{9}, and statistical tests in cluster analysis for non-metric
data \cite{5}.

An interesting topic in the study of random intersection graphs
regards the interrelations between the active and passive graphs.
For example, if $V$ and $W$ represent researchers and research
papers, respectively; and $v\in V$ and $w\in W$ are adjacent if
researcher $v$ is an author of paper $w$. Thus, the resulting active
graph is a collaboration graph on $V$ of researchers, and
accordingly, the passive graph is a relation graph on $W$ of papers.
The correlation of these two graphs shall shed light on the in-depth
architectures and patterns in scientific collaboration \cite{10};
e.g. in what way and to what extent, the appearance of certain
configurations in $G^{active}(n,m,p)$ may influence the topology of
$G^{passive}(n,m,p)$ and vice versa.

In the current paper, as a first step towards this research, we
explore the correlation of vertex degrees between
$G^{active}(n,m,p)$ and $G^{passive}(n,m,p)$, and the corresponding
joint probability generating function is provided by exploiting the
sieve method. Some related works dealing with degree distributions
of random intersection graphs have been done. \cite{1} treats the
degree distribution of $G^{active}(n,m,p)$ with parameters
$m=\lfloor n^{\alpha}\rfloor$ for some $\alpha>0$ and
$p=\sqrt{c}n^{-(1+\alpha)/2}$ for some $c>0$, and the distribution
of the degree of a given vertex is shown to converge to a point mass
at 0, a Poisson distribution or a compound Poisson distribution
depending on whether $\alpha<1$, $\alpha=1$ or $\alpha>1$. \cite{4}
examines the sufficient and necessary conditions for generalized
active random intersection graph to have a Poisson limiting degree
distribution. The related conditions which imply a limiting degree
distribution are given in \cite{6} for generalized passive random
intersection graph. To the best of our knowledge, our work is the
first one devoted to joint degree distributions between the active
and passive graphs. It is obvious that the number of edges in active
and passive graphs should be positively correlated.

The rest of the paper is organized as follows. The main result
(Theorem 1) is presented and proved in Section 2. In Section 3, we
conclude the paper with some further remarks.

\bigskip
\noindent{\Large\textbf{2. Joint probability generating function}}
\smallskip

Let $X$ be the number of vertices of $V\backslash\{v\}$ adjacent in
$G^{active}(n,m,p)$ to a vertex $v\in V$, and $Y$ be the number of
vertices of $W\backslash\{w\}$ adjacent in $G^{passive}(n,m,p)$ to a
vertex $w\in W$. That is, $X$ and $Y$ are typical vertex degrees in
$G^{active}(n,m,p)$ and $G^{passive}(n,m,p)$, respectively. The
joint probability generating function of $X$ and $Y$ is defined to
be $F(x,y):=Ex^Xy^Y$ for $x,y\in\mathbb{R}$. The main result in this
paper is stated as follows.

\medskip
\noindent\textbf{Theorem 1.}\itshape \quad The joint probability
generating function $F(x,y)$ is given by
\begin{eqnarray*}
F(x,y)&=&\sum_{k=0}^{n-1}\sum_{l=0}^{m-1}{n-1 \choose k}{m-1\choose
l}x^{n-1-k}(1-x)^ky^{m-1-l}(1-y)^l\\
& &\cdot
\big[1-p+p(1-p)^k\big]^{m-1-l}\big[1-p+p(1-p)^l\big]^{n-1-k}\\
& &\cdot \Big[(1-p)^{k+l}p+(1-p)\sum_{i=0}^l{l\choose
i}p^{i}(1-p)^{l-i}\big[(1-p)^{i+1}+p(1-p)^l\big]^k\Big].
\end{eqnarray*}
\normalfont\smallskip

Let $F(x)$ and $F(y)$ denote the probability generating functions
for random variables $X$ and $Y$, respectively. Since $F(x)=F(x,1)$
and $F(y)=F(1,y)$, we have the following corollary which is
consistent with Theorem 1 in \cite{1}.

\medskip
\noindent\textbf{Corollary 1.}\itshape \quad The probability
generating functions $F(x)$ and $F(y)$ are given by
$$
F(x)=\sum_{k=0}^{n-1}{n-1 \choose k}x^{n-1-k}(1-x)^k
\big[1-p+p(1-p)^k\big]^{m}
$$
and
$$
F(y)=\sum_{l=0}^{m-1}{m-1 \choose l}y^{m-1-l}(1-y)^l
\big[1-p+p(1-p)^l\big]^{n}
$$
respectively.
 \normalfont\smallskip

To prove Theorem 1, we will invoke a lemma which is a probability
generating function version of the sieve method, whose proof is
similar to Lemma 1 in \cite{1}. In fact, the following lemma can be
viewed as a high dimensional extension of Lemma 1 in \cite{1}. Let
$P_1$ and $P_2$ be two disjoint sets of properties that a random
object can take on. Let $p_{k,l}$ be the probability that the object
takes on exactly $k$ properties in $P_1$ and $l$ properties in
$P_2$. The related probability generating function is defined as
$F(x,y):=\sum_{k,l\ge0}p_{k,l}x^ky^l$.

\medskip
\noindent\textbf{Lemma 1.}\itshape \quad Given $S_1\subset P_1$ and
$S_2\subset P_2$, we define $N_{S_1,S_2}$ to be the event that the
random object possesses properties $S_1$ and $S_2$. Define
$$
N_{k,l}:=\sum_{|S_1|=k,|S_2|=l}P(N_{S_1,S_2})\quad and\quad
N(x,y):=\sum_{k,l\ge0}N_{k,l}x^ky^l.
$$
Then we have
$$
F(x,y)=N(x-1,y-1).
$$
 \normalfont

\smallskip
\noindent\textbf{Proof}. Let $Y_1$ and $Y_2$ be the numbers of
properties that the random object possesses in $P_1$ and $P_2$,
respectively. Let $1_{N_{S_1,S_2}}$ be the indicator function of
event $N_{S_1,S_2}$. We clearly obtain
$$
N_{k,l}=\sum_{|S_1|=k,|S_2|=l}E(1_{N_{S_1,S_2}})=E\bigg(\sum_{|S_1|=k,|S_2|=l}1_{N_{S_1,S_2}}\bigg)=E\bigg({Y_1\choose
k}{Y_2\choose l}\bigg).
$$
Hence,
\begin{eqnarray*}
N(x,y)&=&\sum_{k,l\ge0}E\bigg({Y_1\choose k}{Y_2\choose
l}\bigg)x^ky^l=E\bigg[\sum_{k\ge0}{Y_1\choose
k}x^k\sum_{l\ge0}{Y_2\choose
l}y^l\bigg]\\
&=&E\big[(x+1)^{Y_1}(y+1)^{Y_2}\big]=F(x+1,y+1).
\end{eqnarray*}
$\Box$

In the following proof, we will take $P_1$ as the set of $n-1$
properties consisting of the non-adjacency of a fixed vertex to the
other $n-1$ vertices in $V$; and take $P_2$ as the set of $m-1$
properties regarding $W$ similarly.

\smallskip
\noindent\textbf{Proof of Theorem 1}. For $v\in V$ and $w\in W$, let
$p'_{k,l}$ be the probability that exactly $k$ vertices in
$V\backslash\{v\}$ are not adjacent to $v$ and $l$ vertices in
$W\backslash\{w\}$ are not adjacent to $w$. Let $G(x,y)$ be the
corresponding probability generating function. The probability that
the fixed vertices $v$ and $w$ are adjacent to none of the vertices
represented by $S_1\subset V\backslash\{v\}$ and $S_2\subset
W\backslash\{w\}$ respectively is given by
$$
P(N_{S_1,S_2})=P(N_{S_1,S_2}| v\sim w)\cdot p+P(N_{S_1,S_2}|
v\not\sim w)\cdot(1-p),
$$
where
\begin{eqnarray}
P(N_{S_1,S_2}| v\sim
w)&=&\sum_{i=1}^{n-|S_1|}\sum_{j=1}^{m-|S_2|}{m-1-|S_2|\choose
j-1}{n-1-|S_1|\choose
i-1}\nonumber\\
& &\cdot
p^{j-1}(1-p)^{m-j}(1-p)^{(j-1)|S_1|}p^{i-1}(1-p)^{n-i}(1-p)^{(i-1)|S_2|}\label{1}
\end{eqnarray}
and
\begin{eqnarray}
P(N_{S_1,S_2}| v\not\sim
w)&=&\sum_{i_s=0}^{|S_1|}\sum_{i_o=0}^{n-1-|S_1|}\sum_{j_s=0}^{|S_2|}\sum_{j_o=0}^{m-1-|S_2|}{m-1-|S_2|\choose
j_o}{|S_2|\choose j_s}{n-1-|S_1|\choose i_o}{|S_1|\choose
i_s}\nonumber\\
& &\cdot
p^{j_o+j_s}(1-p)^{m-1-j_o-j_s}p^{i_o+i_s}(1-p)^{n-1-i_o-i_s}\nonumber\\
& &\cdot (1-p)^{(j_o+j_s)|S_1|+(i_o+i_s)|S_2|-i_sj_s}.\label{2}
\end{eqnarray}
In the above expression (\ref{1}), the index $i$ counts the number
of vertices of $V$ adjacent to $w$ in $G^*(n,m,p)$ and similarly,
$j$ counts the number of vertices of $W$ adjacent to $v$ in
$G^*(n,m,p)$; c.f. Fig. 1. For the expression (\ref{2}), $i_o+i_s$
counts the number of vertices of $V$ adjacent to $w$ in
$G^*(n,m,p)$, where $i_o$ counts vertices outside $S_1$ while $i_s$
inside. The roles of indices $j_o$ and $j_s$ can be interpreted
likewise; c.f. Fig. 2. The inclusion and exclusion principle is
utilized here.

Hence we have
\begin{eqnarray*}
N_{k,l}&=&{n-1\choose k}{m-1\choose
l}\bigg[\sum_{i=1}^{n-k}\sum_{j=1}^{m-l}{m-1-l\choose
j-1}{n-1-k\choose
i-1}\\
& &\cdot p^{j-1}(1-p)^{m-j+(j-1)k}p^{i-1}(1-p)^{n-i+(i-1)l}p\\
&
&+\sum_{i_s=0}^{k}\sum_{i_o=0}^{n-1-k}\sum_{j_s=0}^{l}\sum_{j_o=0}^{m-1-l}{m-1-l\choose
j_o}{l\choose j_s}{n-1-k\choose i_o}{k\choose i_s}\\
& &\cdot
p^{j_o+j_s}(1-p)^{m-1-j_o-j_s+(i_o+i_s)k}p^{i_o+i_s}(1-p)^{n-1-i_o-i_s+(i_o+i_s)l-i_sj_s}
(1-p)\bigg]\\
&=&{n-1\choose k}{m-1\choose
l}\bigg[\big[1-p+p(1-p)^l\big]^{n-1-k}(1-p)^k\big[1-p+p(1-p)^k\big]^
{m-1-l}\\
& &\cdot(1-p)^lp+\big[1-p+p(1-p)^k\big]^{m-1-l}\big[1-p+p(1-p)^l\big]^{n-1-k}(1-p)\\
& &\cdot \sum_{i_s=0}^{k}\sum_{j_s=0}^l{k\choose i_s}{l\choose
j_s}\big[p(1-p)^l\big]^{i_s}(1-p)^{k-i_s}\big[p(1-p)^k\big]^{j_s}(1-p)^{l-j_s}(1-p)^{-i_sj_s}\bigg]\\
&=&{n-1\choose k}{m-1\choose
l}\big[1-p+p(1-p)^k\big]^{m-1-l}\big[1-p+p(1-p)^l\big]^
{n-1-k}\\
& &\cdot\bigg[(1-p)^{k+l}p+(1-p)\sum_{i=0}^{l}{l\choose
i}p^i(1-p)^{l-i}\big[(1-p)^{i+1}+p(1-p)^l\big]^k\bigg].
\end{eqnarray*}
Consequently, Lemma 1 yields $G(x,y)=N(x-1,y-1)$, where
\begin{eqnarray*}
N(x,y)&=&\sum_{k=0}^{n-1}\sum_{l=0}^{m-1}N_{k,l}x^ky^l\\
&=&\sum_{k=0}^{n-1}\sum_{l=0}^{m-1}x^ky^l{n-1\choose k}{m-1\choose
l}\big[1-p+p(1-p)^k\big]^{m-1-l}\big[1-p+p(1-p)^l\big]^
{n-1-k}\\
& &\cdot\bigg[(1-p)^{k+l}p+(1-p)\sum_{i=0}^{l}{l\choose
i}p^i(1-p)^{l-i}\big[(1-p)^{i+1}+p(1-p)^l\big]^k\bigg].
\end{eqnarray*}
The result then follows from the fact that
$$
F(x,y)=\sum_{k=0}^{n-1}\sum_{l=0}^{m-1}p_{k,l}x^ky^l=\sum_{k=0}^{n-1}\sum_{l=0}^{m-1}p'_{k,l}x^{n-1-k}y^{m-1-l}=x^{n-1}y^{m-1}G(x^{-1},y^{-1}).
$$
$\Box$

\bigskip
\noindent{\Large\textbf{3. Concluding remarks}}
\smallskip

In this letter, we study the joint probability generating function
$F(x,y)$ for degrees of $G^{active}(n,m,p)$ and $G^{passive}(n,m,p)$
by employing the sieve method. We mention that our work is only a
preliminary step to the interesting and meaningful topic of
correlation of active/passive random intersection graphs. Corollary
1 gives the probability generating functions of the marginal
distributions $P(X=k)$ and $P(Y=l)$, respectively. As is known, the
random variables $X$ and $Y$ are independent if and only if
$F(x,y)=F(x)F(y)$ for all $x,y$ (c.f. \cite{7} pp. 279). Hence, we
may tackle the interdependence between the active graph and passive
graph through their probability generating functions. For the future
research, the method developed in this paper may be applied to the
joint distribution of degrees of generalized active and passive
graphs in \cite{5}. A more hard question is to investigate the joint
distribution of clusters in these two models.

\bigskip

\begin{figure}[!t]
\centering \scalebox{0.5}{\includegraphics{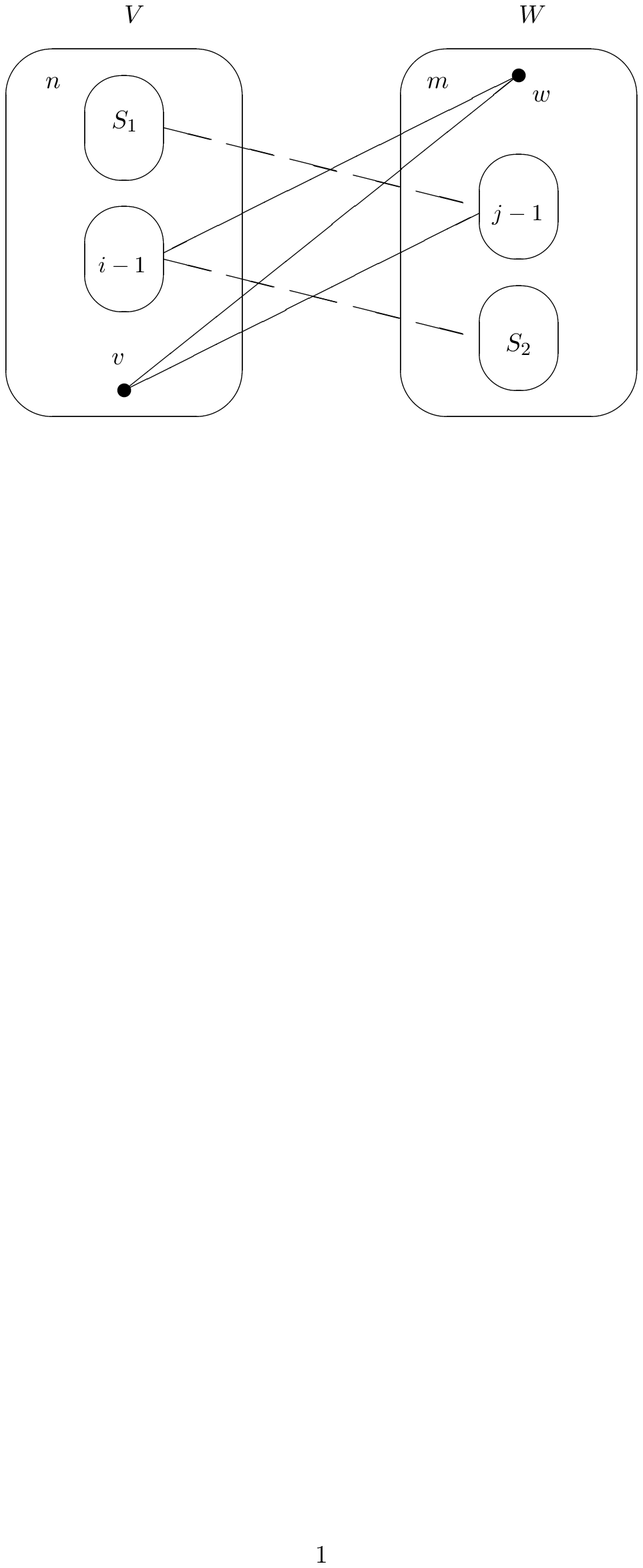}}
\caption{$G^*(n,m,p)$ in the case of $v\sim w$. Solid line
represents edges and dashed line represents non-edges.}
\label{fig_sim}
\end{figure}

\begin{figure}[!t]
\centering \scalebox{0.5}{\includegraphics{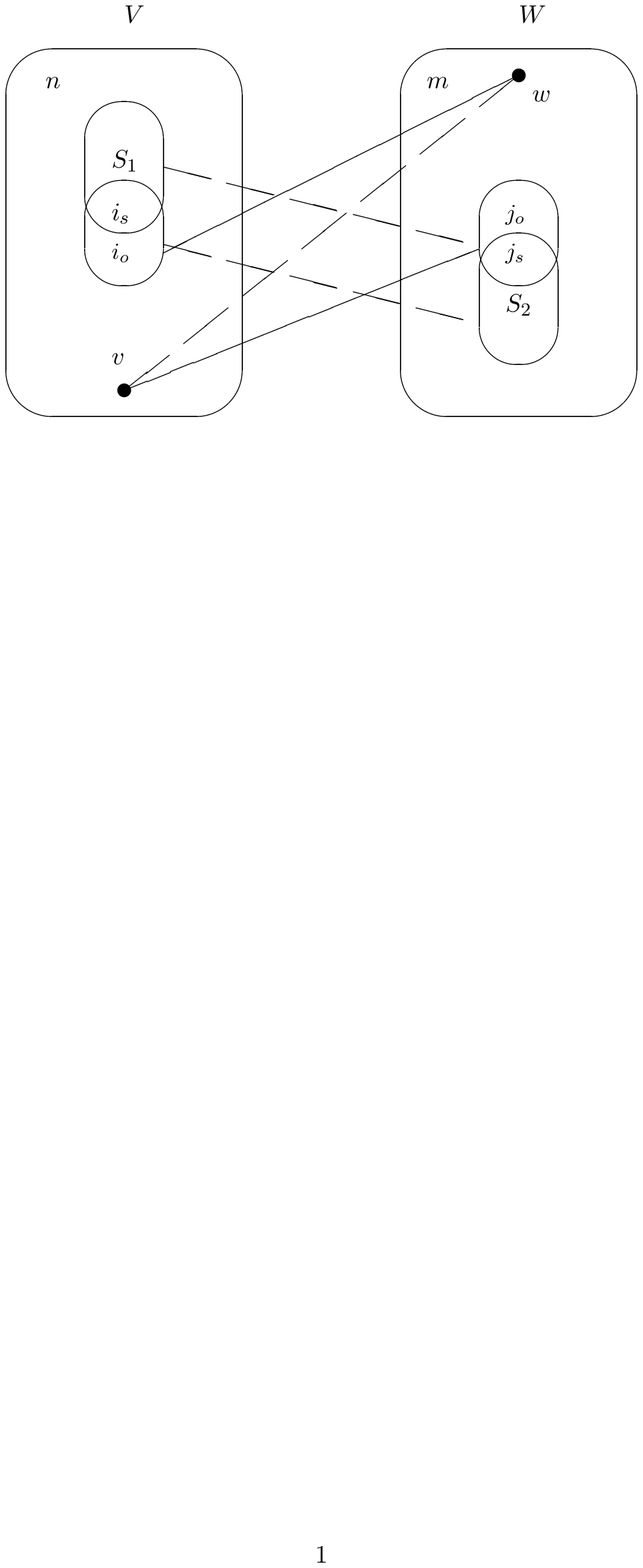}}
\caption{$G^*(n,m,p)$ in the case of $v\not\sim w$. Solid line
represents edges and dashed line represents non-edges}
\label{fig_sim}
\end{figure}

\end{document}